\documentclass[journal]{IEEEtran}
\usepackage{graphicx}
\usepackage{subcaption}
\ifCLASSINFOpdf
\else
\fi
\usepackage{xcolor}
\usepackage{amsmath}
\usepackage{booktabs}
\usepackage{multirow}
\usepackage{hyperref}
\usepackage{algorithm} 
\usepackage{cite}
\usepackage{algpseudocode}

\begin{document}
%


\title{Modeling to Generate Alternatives for Robustness of Mixed Integer DC Optimal Power Flow}

%
%
%

\author{Constance Crozier,~\IEEEmembership{Member,~IEEE}
\thanks{C. Crozier is with the H. Milton School of Industrial and Systems Engineering, Georgia Institute of Technology, Atlanta, GA, 30332 USA e-mail: ccrozier8@gatech.edu.}}

\maketitle

\begin{abstract}
Transmission system operators face a variety of discrete operational decisions, such as switching of branches and/or devices. Incorporating these decisions into optimal power flow (OPF) results in mixed-integer non-linear programming problems (MINLPs), which can't presently be solved at scale in the required time. Various linearizations of the OPF exist, most famously the DC-OPF, which can be leveraged to find integer decisions. However, these linearizations can yield very poor integer solutions in some edge cases, making them challenging to incorporate into control rooms. This paper introduces the use of modeling to generate alternatives (MGA) to find alternative solutions to the linearized problems, reducing the chance of finding no AC feasible solutions. We test this approach using 13 networks where the DC linearization results in infeasible integer decisions, and MGA finds a solution in all cases. The MGA search criteria selected drastically affects the number and quality of solutions found, so network specific search functions may be necessary. 

\end{abstract}

\begin{IEEEkeywords}
AC Power flow, Modelling to generate alternatives, Transmission switching, Unit commitment
\end{IEEEkeywords}

%
\IEEEpeerreviewmaketitle

\section*{Nomenclature}
\textbf{Decision variables}\\[.1cm]
\begin{tabular}{ll}
$p_g$ & Device power consumption/generation (MW)\\
$p_e$ & Line power flow (MW)\\
$\theta_i$ & Voltage angle at bus $i$ (radians)\\
$x_e$ & On/off status of branch $e$\\
$x_g$ & On/off status of device $g$\\
\end{tabular}

\vspace{.1cm}

\textbf{Parameters}\\[.1cm]
\begin{tabular}{ll}
$c_g$ & Cost of device $g$ (\$/MW)\\
$b_e$ & Susceptance of branch $e$\\
$M$ & Large number\\
$n_{on}$ & Number of connected branches\\
$n_{off}$ & Number of disconnected branches\\
$\beta_e$ & (Random) weighting applied to branch $e$\\
$\delta_f$ & Tolerated increase in objective function (\$)\\
\end{tabular}

\vspace{.1cm}

\textbf{Sets}\\[.1cm]
\begin{tabular}{ll}
$E_{io}$ & Branches with origin connected to bus $i$\\
$E_{id}$ & Branches with destination connected to bus $i$\\
$E_{off}$ & Branches which are turned off\\
$E_{on}$ & Branches which are turned on\\
$G_i$ & Devices connected to bus $i$\\
$G_{off}$ & Devices which are turned on\\
$G_{on}$ & Devices which are turned on\\
\end{tabular}

\vspace{.1cm}

\textbf{Notation}\\[.1cm]
\begin{tabular}{ll}
$e_o$ & Bus origin of branch $e$\\
$e_d$ & Bus destination of branch $e$\\
\end{tabular}

\begin{tabular}{ll}
$^{max}$ & Maximum value\\
$^{min}$ &Minimum value\\
$^*$ &At the current optimal solution\\
\end{tabular}


\section{Introduction}
%
%
%
%
\IEEEPARstart{T}{his} paper investigates the use of modeling to generate alternatives to increase the robustness of linearized methods for mixed integer optimal power flow formulations. 

Power system operators are concerned with maintaining a reliable, low cost supply of electric power. This involves solving a set of problems referred to as optimal power flow (OPF), which seek the minimum cost system subject to a variety of physics and operational constraints. In practice many operators rely on the linearized direct current (DC) approximation, rather than the non-linear alternating current (AC) equivalent. While the classic problems studied are largely continuous, operators face a number of discrete decisions such as switching lines, committing generation units, transformer tap configurations, and others. 

The majority of the research investigating discrete actions utilizes the DC approximation, resulting in tractable mixed integer linear programming (MILP) problems~\cite{10080945,4492805,4957010,4558426}. However, the DC linearization makes several approximations which mean that the solutions are non-physical, and necessary do not meet the (more realistic) AC power flow constraints~\cite{baker2021solutions}. Furthermore, several recent works have demonstrated that the DC approximation sometimes lead to poor integer decisions, which produce no feasible solution~\cite{6646289,7038446}. Alternative approaches attempt to solve the AC mixed integer non-linear programming (MINLP) formulations, utilizing decomposition~\cite{6253283,LAN2021107140}, relaxations~\cite{7849135}, or bi-level formulations~\cite{7494647}. However these have so far only shown success on small simplified cases. Heuristic search methods have also shown promise~\cite{7038445,6758897,CROZIER2022100628}, although these typically involve solving a very large number of AC OPF problems, so can be prohibitively slow and expensive.

Many of the limitations of the DC MILP formulations are based on edge cases; for many networks and loading scenarios they produce good results, but in some cases the integer decisions are infeasible. In this paper, we instead propose using a method that purposefully explores the DC MILP solution space, thus increasing the chance of finding an optimal AC solution.

Modeling to generate alternatives (MGA) is a technique for generating a small number of different solutions for decision makers to consider~\cite{105076}. In this method a secondary optimization problem is constructed which, alongside the original constraints, requires the original objective function is within a tolerance of the optimal. New objective functions are then explored to identify different solutions; for example hop, skip, and jump explicitly searches for maximally different solutions, while the random vector approach uses random linear weights in the objective -- both of these are described in greater detail later in Section III. 

Originally, MGA was proposed for zoning problems such as the forestry. However, in recent years, it has been applied to energy sector planning problems~\cite{DECAROLIS2011145,lau2024measuring}. These works have mostly focused on sizing of generation, storage, and transmission -- for example in Europe~\cite{NEUMANN2021106690}, the United States~\cite{DECAROLIS2016300}, and for global energy pathways~\cite{PRICE2017356}. Energy system planning is a well-suited to MGA because there can be multiple systems which have similar projected costs but very different technology mixes. While traditional optimization approaches will yield only a single solution, there may be reasons that one solution may be preferable over the other (e.g. including more nuclear for political favor and security-of-supply). To the author's knowledge, MGA has been exclusively used for planning problems in the energy sector. 

In this paper we instead explore the use of MGA for transmission system operation. 
By using MGA to find multiple similar, but distinct, solutions in the linearized space, we hope to increase the quality of solutions obtained. Transmission system operation also often involves some degree of operator-specific expertise, so presenting multiple possible solutions to operators may be further desirable. The author therefore believes that there will be further opportunities for MGA to improve transmission system operation, which extend beyond the specific problems considered here.


The contributions of this paper can be summarized as follows. Firstly, that we develop an alternative use case for MGA -- where alternative solutions to the approximate problem are explored to improve chances of finding a globally optimal solution. Second, that we demonstrate this paradigm on both the optimal transmission switching and unit commitment problems, potentially allowing a suite of previously developed linearized methods to be practically realizable. We demonstrate that this approach can lead to a high quality solution much faster than popular heuristic methods. Finally, that we explore the efficacy of various MGA search criteria on a variety of different networks.


The remainder of the paper is structured as follows. In Sections \ref{sec:ots} and \ref{sec:mga} we will introduce the optimal power flow and MGA formulations respectively. Section \ref{sec:prop} details the proposed implementation of MGA for solution exploration in the OTS problem. Section \ref{sec:test} details the testing protocol, used to create the results in Section \ref{sec:res}, and Section \ref{sec:con} concludes the paper. 





\section{Mixed Integer Optimal Power Flow Problems} \label{sec:ots}

In this section, we describe both the AC and DC optimal power flow problems in their continuous form. Then we describe example integer decisions faced by power system operators, and discuss why MILP formulations may not produce good results. 

\subsection{DC Optimal Power Flow}

The classical DC optimal power flow problem (DC-OPF) aims to find the minimum cost operation of a power transmission network and can be expressed as:

\begin{subequations}
\begin{align}
\min_{\theta, p_g}\quad& \sum_g c_g p_g \label{eq:obj}\\
\text{s.t.}\quad & 0 = \sum_{g\in G_i} p_g - \sum_{i\in E_{io}} p_e- \sum_{i\in E_{id}} p_e \quad \forall i \label{eq:pb}\\
& p_e = b_e (\theta_{e_o} - \theta_{e_d})\quad \forall e \label{eq:dc}\\
& -p_e^{max} \leq p_e \leq p_e^{max} \quad \forall e \label{eq:e_bound}\\
& p_g^{min} \leq p_g \leq p_g^{max} \quad \forall g \label{eq:g_bound}
\end{align}
\end{subequations}

The objective function \eqref{eq:obj} describes the total costs from generation of power; $p_g$ represent piece-wise linear segments of generation or load which each come at a constant cost of $c_g$ \$/MWh. Note that the traditional formulation includes only controllable generation, however here we take a more modern approach which includes price sensitive demand. In practice, we can consider three types of devices: generators, dispatchable, and non-dispatchable loads. For generators $p_g$ will be positive, for loads $p_g$ will be negative and if the loads are non-dispatchable $p_g^{min}=p_g^{max}$ such that their value is fixed. 

The constraint \eqref{eq:pb} enforces Kirchoff's circuit law, that the net current (here power) into each node should be zero. The power flowing down branch $e$, $p_e$, is not directly controllable, but dictated by the physics of the circuit. In this case we use the DC approximation \eqref{eq:dc} which expresses the power as linear to the difference in voltage angle between the origin and destination buses. Finally we enforce generator/load bounds \eqref{eq:g_bound} as set by the devices, and branch flow limits \eqref{eq:e_bound} which avoid overheating. Although the only controllable variable is device output, $p_g$, we also take the bus voltage angles, $\theta$, as decision variables to avoid complex non-linear constraints.

\subsection{AC Power Flow Constraints}

The DC power constraints make several simplifying assumptions, notably: that all voltage magnitudes are unity, ignoring effects of reactive power, and utilizing the small angle approximation. A more accurate expression for the $p_e$ is found using the AC power flow constraints can be written as follows:
\begin{align}
p_e &= v_{e_o}v_{e_d} \Big ( b_e \text{sin}(\theta_{e_o} - \theta_{e_d}) +  g_e \text{cos}(\theta_{e_o} - \theta_{e_d})\Big) \label{eq:ac_pe}
\end{align}
where $v$ give the voltage magnitudes, and $g_e,b_e$ define the electrical parameters of the line or transformer $e$. This expression adds significant complexity to the OPF problem; given that $v,\theta$ are both decision variables, we end up with a quartic multiplication of decision variables as well as sinusoidal operators. Note that under AC, separate expressions are also required for the power flowing into and out of a branch due to losses. 

Additionally, we include expressions for the balance of reactive power, and the reactive power flowing down a branch:
\begin{align}
0 &= \sum_{g\in G_i} q_g - \sum_{i\in E_{io}} q_e- \sum_{i\in E_{id}} q_e \quad \forall i \label{eq:qb}\\
q_e &= v_{e_o}v_{e_d} \Big ( g_e \text{sin}(\theta_{e_o} - \theta_{e_d}) -  b_e \text{cos}(\theta_{e_o} - \theta_{e_d})\Big) \label{eq:ac_qe}\,,
\end{align}
where $q_g$ are the (relatively small) reactive power generation or consumption of devices -- these are bounded decision variables. The added complexity makes it significantly more challenging to solve the AC-OPF problem, a mixed-integer non-linear programming problem. Current computational tools can not solve these problems at realistic scales.

A commonly proposed approach is to optimize the topology using the linearized model, and then to run an AC-OPF with the chosen topology. This approach is later referred to as an AC recovery problem, as it uses changes in generation and load to offset the error from the DC constraints. However, as we will show later, the are cases where a chosen topology has no feasible solution using the AC constraints.

\subsection{Common Mixed Integer Problems}

The formulations described above are continuous, with the decision variables all reflecting quantities which can take a smooth range of values. In practice, transmission system operators face a number of integer decisions which also need to be incorporated into the optimal power flow problem. Here we describe some common problems, and discuss why the DC linearization does not always provide good integer decisions.

\subsubsection{Transmission Switching}

Optimal Transmission Switching (OTS) extends the OPF problem to include switching on/off of lines or transformers. These switches are routinely used to restore network redundancy after a component outage, ease congestion, or manage power quality. The DC-OTS problem can be expressed as: 
\begin{subequations}
\begin{align}
\min_{\theta, p_g, x_e}\quad& \sum_g c_g p_g\\
\text{s.t.} \quad&\eqref{eq:pb},\eqref{eq:e_bound},\eqref{eq:g_bound}\\
& p_e \leq b_e (\theta_{e_o} - \theta_{e_d}) + M(1-x_e)\quad \forall e \label{eq:bigm1}\\
& p_e \geq b_e (\theta_{e_o} - \theta_{e_d}) - M(1-x_e)\quad \forall e \label{eq:bigm2}\\
& -M x_e\quad \leq p_e \leq M x_e\quad \forall e \label{eq:bigm3}\\
& \sum_{e \in E_{io}, E_{id}} x_e \geq \text{min} (2, \sum_{e \in E_{io}, E_{id}} 1) \quad \forall i \label{eq:islands} \\
& \theta^{min} \leq \theta_e \leq \theta^{max} \quad\forall e
\end{align}
\end{subequations}  

The objective function remains unchanged, but we have introduced new variables $x_e$ which are binary variables describing whether branch $e$ is on or off. In addition to the constraints from the DC-OPF formulation, we include constraints to enforce the relationship between a branch's on/off position and the power flowing down it; \eqref{eq:bigm1}, \eqref{eq:bigm2} and \eqref{eq:bigm3} enforces the logical statement:
\begin{align*}
p_e &= b_e (\theta_{e_o} - \theta_{e_d}) \quad \text{if  }x_e = 1\\
&= 0 \quad\quad \quad \quad \quad \quad \text{otherwise}
\end{align*}
The constraints \eqref{eq:bigm1}, \eqref{eq:bigm2} exploit the common big-$M$ formulation, where $M$ is a relatively large number. Where $x_e=1$ the $M$ term disappears and we are left with two opposing inequalities (thus enforcing equality); where $x_e=0$ the constraints become very loose (necessarily inactive assuming $M>p_e^{max}$). Similarly, \eqref{eq:bigm3} will enforce $p_e=0$ if $x_e = 0$ and become inactive if $x_e=1$. 

A common desire when we consider allowing disconnection of transmission lines is to maintain connectivity of the network, i.e. to avoid creating islands (splitting the network into two unconnected grids). In this formulation we include a heuristic constraint to maintain good connectivity \eqref{eq:islands}, which ensures that each bus must be connected to at least two others where possible. This constraint is a heuristic, in that there are both examples of islands satisfy the expression, and connected networks which violate. Other heuristics are available, and in many applications an anti-islanding constraint is not necessary in practice. However, for the cases considered in this paper, \eqref{eq:islands} was necessary to avoid islands and did so successfully.

Note that in practice, system operators may also place additional constraints on the allowable switches -- for example, an upper bound on the total number of switches from the network position using:
\begin{align}
    \sum_e (x_e - x_e^{0}) \leq n_{s}^{max}\,,
\end{align}
where $x_e^{0}$ is the initial position of branch $e$ and $n_{s}^{max}$ is the maximum allowable switches. Alternatively, they may wish to prevent particular combinations of switches using:
\begin{align}
    x_{e_1} + x_{e_2} \leq 1\,,
\end{align}
which ensures that a maximum of 1 of branches $e_1$ and $e_2$ are switched. A variety of other specific constraints are possible, but all can be expressed linearly in $x$ using logic.


\subsubsection{Unit Commitment}

In the classic OPF problem we assume that generators can take a continuous range of real power outputs, from $p_g^{min}$ to $p_g^{max}$. Often, to ensure feasibility, it is assumed that $p_g^{min}=0$. However, in practice generators have a non-zero lower bound for the power they can output while maintaining reasonable quality, below this the generator must be shut down. The unit commitment (UC) problem considers which generators should be turned on, normally over multiple periods. The single period DC-UC problem can be expressed as:
\begin{subequations}
\begin{align}
\min_{\theta, p_g, x_e}\quad& \sum_g c_g p_g\\
\text{s.t.} \quad&\eqref{eq:pb},\eqref{eq:e_bound},\eqref{eq:dc},\eqref{eq:g_bound}\\
& p_g \leq p_g^{max}x_g\quad \forall g \label{eq:newbound}\\
& p_g \geq p_g^{min}x_g\quad \forall g \label{eq:newbound}\,,
\end{align}
\end{subequations} 
where $x_g$ describes whether each generator is on or off. These decisions become more important when the problem is extended to its multi-period case, due to ramping limits e.g.:
\begin{align}
    p_{g_t} &\leq p_{g_{t-1}} + \delta_g\\
    p_{g_t} &\geq p_{g_{t-1}} - \delta_g\,
\end{align}
where $p_{g_t}$ describes the value of $p_g$ at time $t$ and $\delta_g$ is the generator specific maximum ramp in a single timestep. The inclusion of ramp rates and multi-period planning significantly increases the complexity of unit commitment decisions, especially in systems where the ramp rates are low and the demand varies significantly.

\subsection{Limitations of Mixed-Integer DC Decisions}

The DC power flow constraints are a linear approximation of the AC non-linear equivalent. More complicated linearizations are available, but they tend to fit better over small operating regions but worse over larger load changes~\cite{9916903}. The solutions to linear optimization problems lie at vertices of the feasible space, or intersection of constraints. This means that for DC-OPF problems the solutions typically result in two or more of the constraints being active: typically either thermal limits and/or generator upper bounds. However, the DC thermal limits are an approximation the solution often includes a set of generator outputs which are not feasible in AC. For continuous problems, it is typically easy to apply corrections to the DC-OPF solution such that it conforms to the AC physics. However, for DC-MILP formulations the integer decisions are typically fixed and can (in some cases) result in a zero feasible set. 

\begin{figure}
    \centering
    \includegraphics[width=0.7\columnwidth,trim={1cm 0.9cm 1cm 1cm},clip]{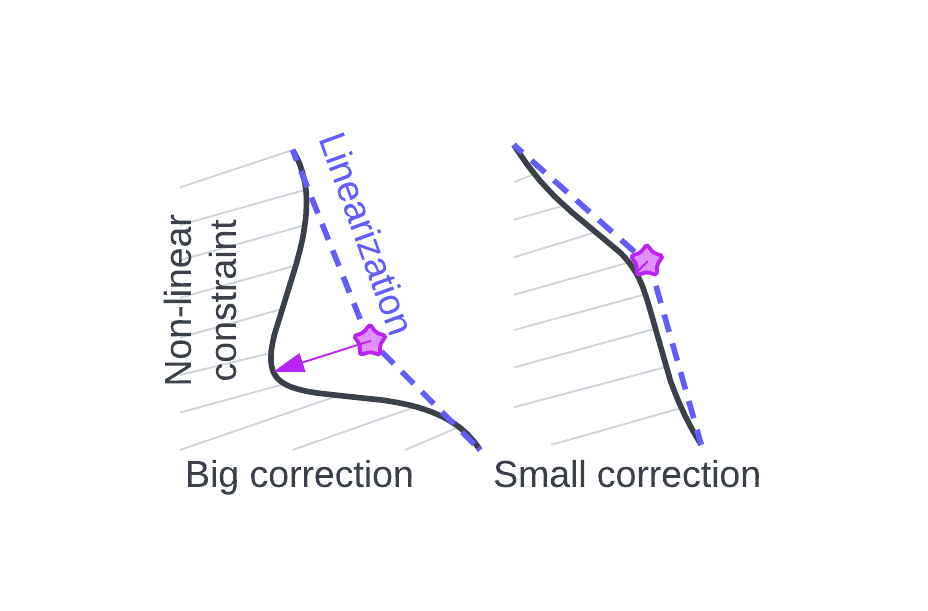}
    \caption{Solutions to the linearized problem lie on vertices, the distance from the feasible non-linear region will vary.}
    \label{fig:lin}
\end{figure}

In practice, researchers and professionals find that for some networks and loading scenarios the DC constraints fit well, and the MILP formulations return good results. This is because in some cases the optimal solution to the DC-MILP exists at the intersection of constraints which closely represent the true feasible space. This effect is visualized in Fig. \ref{fig:lin}. 

\section{Modeling to Generate Alternatives} \label{sec:mga}

Traditional optimization approaches converge towards a single solution. This paradigm relies on the problem having a deterministic objective function, and well defined constraints. In practice, objective functions typically contain some degree of uncertainty regarding cost parameters. Furthermore, there may be some solutions which are preferred over others for reasons that are challenging to define mathematically (e.g. political favor). There may be a large number of feasible solutions to the optimization problem which are \textit{nearly-optimal}, whose difference in objective is considerably smaller than the uncertainty in problem definition. 

Modeling to generate alternatives (MGA) is a method of explicitly searching for alternative, near-optimal solutions. This allows a number of options to be presented to a human decision-maker, who may value some solutions over others for reasons that are hard to express mathematically. 
There are two traditional methods for searching for alternative solutions: hop, skip, and jump (HSJ) and random vector. In this section we will define these algorithms in reference to the DC-OTS problems, as well as proposing two alternative algorithms that may be better suited to the OTS problem.

\subsection{Hop, Skip, and Jump}

The most classical MGA search method is referred to as hop, skip, and jump (HSJ). In the first iteration of HSJ applied to the DC-OTS problem can be written as:
\begin{subequations}\label{eq:mga}
\begin{align}
\min_{\theta, p_g, x_e}\quad& f_0 = \sum_{e \in E_{on}^*} x_e \label{eq:mgaf}\\
\text{s.t.}\quad& \sum_g c_g p_g \leq \sum_g c_g p_g^* + \delta_f \label{eq:qual}\\
& \eqref{eq:pb},\eqref{eq:e_bound},\eqref{eq:g_bound},\eqref{eq:bigm1},\eqref{eq:bigm2},\eqref{eq:bigm3}, \eqref{eq:islands}
\end{align}
\end{subequations}  
where $E_{on}^*$ is the set of connected branches in the DC-OTS topology, and $p_g^*$ is the generator/load outputs at the DC-OTS solution. The objective function sums the non-zero integer variables from the previous solution. The intuition behind this is that it incentivizes solutions not using those branches. The quality of the solution is maintained by \eqref{eq:qual}, which ensures the transmission costs must be within $\delta_f$ of the optimal value found in the original solution. The size of $\delta_f$ will affect the solution search space, with larger values offering greater exploration but potentially worse solutions. For this application it is likely that a value of $\delta_f=0$ would still yield a large number of potential results. This is because there are typically large number of topologies that would support the same mix of generators and demand. The remainder of the problem constraints match those from the original DC-OTS problem. For the DC-UC problem, the objective function \eqref{eq:mgaf} would change to: $\sum_{g \in G_{on}^*} x_g$ where $G_{on}$ is the set of generators on in the initial solution.

Once the first iteration of HSJ has been completed, the problem objective is updated to include the non-zero elements of the new solution, as:
\begin{align}\label{eq:f_update}
    f_{n} &= f_{n-1} + \sum_{e \in E_{on}^{*n-1}}x_e\,,
\end{align}
where $f_{n}$ represents the objective in the $n^{\text{th}}$ iteration of HSJ and $E_{on}^{*n-1}$ is the set of branches on in the optimal solution to the $n-1^{\text{th}}$ iteration. Crucially, the previous elements remain (otherwise our next solution would likely be the original one again). Iterations continue until a non-unique solution is obtained. The search algorithm can therefore be summarized as follows:
\begin{enumerate}
    \item Generate initial solution of the original problem, $x_0$.
    \item Formulate the initial MGA problem as in \eqref{eq:mga} and generate the next solution, $x_1$.
    \item Update the MGA objective $f_n$ using \eqref{eq:f_update} and resolve the problem; repeat until either a repeated solution is generated, or the maximum iterations are exceeded.
\end{enumerate}

\subsection{Random Vectors}

Alternatively, the random vectors approach seeks to better explore the feasible space, using randomly generated weights:
\begin{align}
f_{n} = \quad& \sum_{e} \beta_e x_e
\end{align}
In this approach each iteration of MGA is independent, meaning that the problem easily lends itself to parallel computing. However, there is no convergence criteria and we do not explicitly seek to find very different solutions, so sufficient exploration of the feasible space can take time. In this paper we implement random vectors using a pre-defined number of iterations (specifically 30). However, to maximize the difference the random weights are generated using latin-hypercube sampling. This is a sampling technique which, given a fixed number of samples, seeks to cover as wide an area of multi-dimensional space as possible~\cite{loh1996latin}.

\subsection{HSJ with Negative Weights}

Unlike previous applications of MGA, there is no summation constraint of the binary variables in the transmission switching problem. This means that (as demonstrated later) HSJ is likely to reward a less connected power networks. This is not always desirable, as operators typically prefer redundancy in network paths and less connected networks are (anecdotally) more prone to infeasible AC problems. 

Therefore, here we first consider the reverse strategy: to reward turning branches on which are not used in the previous solution. This would have the alternative objective:
\begin{align}
f_n = f_{n-1} + \sum_{e \in E_{off}^*}  -x_e\,,
\end{align}
and we would use the same iterative procedure and convergence criteria as with standard HSJ.

\subsection{HSJ with Zero Bias}

In order to maximize exploration, it may be better to not incentivize either more or fewer total connected lines; either option will lead to decreased diversity of solutions. Instead we propose a new modification which explicitly doesn't incentivize a net change in number of lines. We achieve this by ensuring that each added function would return a value of zero for a fully connected network. The objective function update takes the form:


\begin{align}
f_n = f_{n-1} + \sum_{e \in E_{on}^*} \frac{1}{n_{on}^*}x_e - \sum_{e \in E_{off}^*} \frac{1}{n_{off}^*}x_e
\end{align}

\section{Proposed use of MGA for solution exploration} \label{sec:prop}

\begin{figure*}[h!]
    \centering
    \includegraphics[width=18cm,trim={0cm 5cm 0 3.5cm},clip]{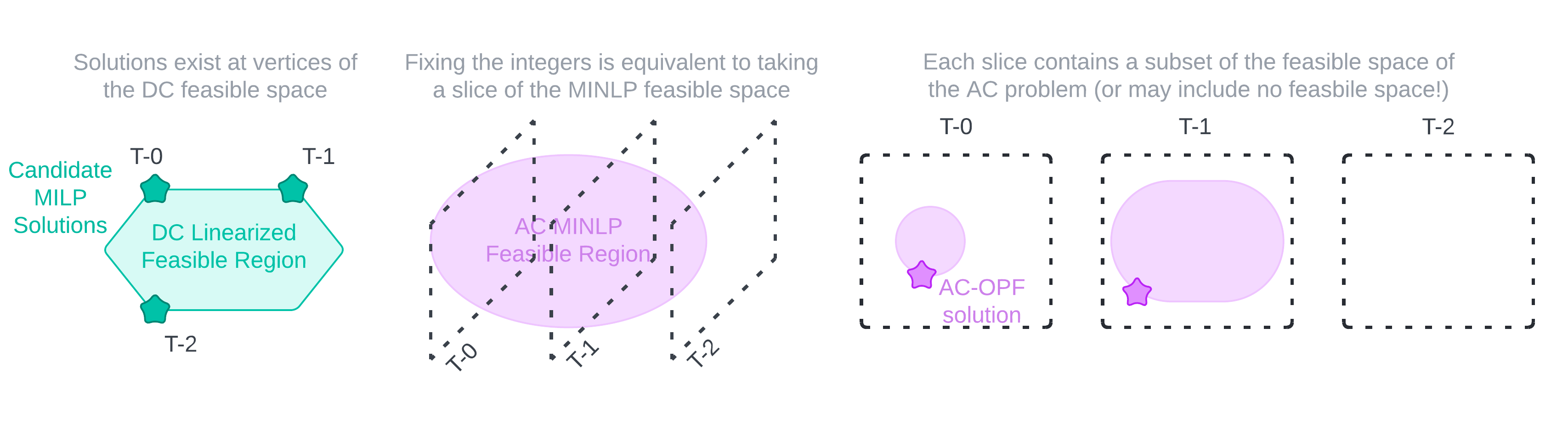}
    \caption{Solutions to the linearized DC-OTS problem exist at the vertices of the feasible polygon. However the AC-OTS problem is higher dimension and bounded by non-linear constraints. Fixing a topology can be thought of as taking a slice of the AC-OTS space. Taking multiple slices increases the chance of finding a good solution to the AC-OTS problem.}
    \label{fig:diag}
\end{figure*}

Traditional MGA applications focus on finding different solutions with similar objective values. In this application, we propose a different paradigm; the use of MGA to improve the chances of finding a good solution using an approximate (linearized) problem. 

The infeasibility of DC-OPF solutions stem from multiple physical processes, including: (1) the effects of losses are ignored, and  (2) voltage limits are not considered. The first approximation means that the DC solution necessarily underestimates both generation required and the branch flows, but these are typically $<5\%$ errors and can often be solved by adjusting generation. The second approximation can result in much larger errors; essentially once the voltages are considered there may be significantly less power that can be safely routed through a section of the network. In these cases, changes in generation may be unable to solve the infeasibility (especially if ramping constraints are included, which we have neglected in this paper). However, the voltage limits only affect some solutions, in \cite{crozier2022data} it is shown that the physical properties of the network mean that only a small number of nodes can experience voltage issues. This means that some linearized solutions likely provides a reasonable starting point.

Here we propose the use of MGA for identifying many possible solutions to the DC-OTS problem with the same, or similar, objective functions. The logic behind this is that the larger the number of candidate solutions, the more likely that one of them is not significantly affected by voltage bounds. The proposed process is visualized in Fig. \ref{fig:diag}. In the first step MGA is performed on the linearized DC-OTS problem, identifying several candidate solutions. Next, for each solution we enforce the chosen topology in an AC-OPF problem, this can be thought of as taking cutting-planes in the higher-dimensional AC-OTS problem. Finally, the AC-OPF problems are run to find a local optima within the plane -- which in some cases may be have a null feasible set.

The AC-OPF problems run can be thought of as a recovery problem -- trying to find a feasible solution from an infeasible starting point. For example, for the OTS problem this can be expressed as follows: 

\begin{subequations}
\begin{align}
    \min_{\theta,v,p_g} \quad &\sum_g c_g p_g\\
    \text{s.t.}\quad &\eqref{eq:pb},\eqref{eq:e_bound},\eqref{eq:g_bound},\eqref{eq:qb}\\
    &\eqref{eq:ac_pe},\eqref{eq:ac_qe} \quad \forall e \in E_{on}^*\\
    &v_{min} \leq v_i \leq v_{max}
\end{align}
\end{subequations}

Note that the line switch decision variables $x_e$ has been removed from the problem, but voltage variables $v$ have been added. Instead of the DC branch flow expression, we impose the AC one \eqref{eq:ac_pe} for all branches that are connected in the chosen topology. Crucially, we allow the generator and loads to vary to the full extent of the original solution; this means that the generation and loads may deviate far from their optimized values, but that if there is a feasible solution using the chosen topology it will be found. In practice, given that there many not be a feasible solution, it is better to made at least the thermal constraints soft -- penalized in the objective function. This can be achieved using additional variables $\delta_e\geq0$ to describe any violation of the thermal constraint, such as:
\begin{align}
    -\delta_e -p_e^{max} \leq p_e \leq p_e^{max}+\delta_e\
\end{align}
These variables are then included in the objective with large positive weights such that the objective becomes $\sum+g c_gp_g + M\sum_e\delta_e$, where $M\gg c_g$. This means that the solver which focus on reducing violations to zero, and once a $\delta_e=0$ solution is achieved, the economic objective is considered. 

\section{Simulation \& Testing Platforms} \label{sec:test}

\begin{figure}[h!]
    \centering
    \includegraphics[width=\linewidth,clip,trim={0 .2cm 0 .5cm}]{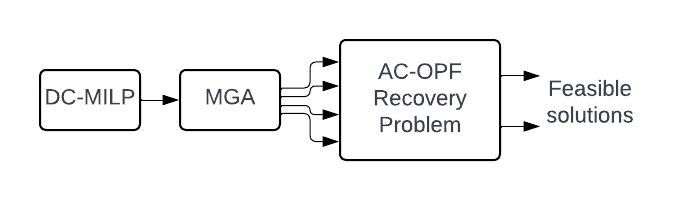}
    \caption{Visualization of the process to find AC feasible solutions}
    \label{fig:flow}
\end{figure}

 In this section, we describe the simulation and testing platform used to assess the performance of MGA for finding good AC integer solutions -- summarized in Fig. \ref{fig:flow}.

We implemented the DC-MILP formulations in Python using cvxpy~\cite{cvxpy} with Gurobi~\cite{gurobi}. Candidate solutions generated by the DC-MILP and MGA iterations are then passed to the AC optimal power flow (AC-OPF) recovery problem. This was implemented in Cython directly interfacing with the IPOPT solver~\cite{ipopt}. Several, acceleration techniques were deployed including analytic forms of the Jacobian and Hessian -- more details are in~\cite{sharadga2024optimizing}.

For testing we utilized networks from the ARPA-E Grid Optimization competition~\cite{go_networks}. This competition aimed to bridge the gap between industry and academic research on power grid optimization; the datasets released represent realistic (but not real) networks. Specifically we investigated scenarios from two networks that were flagged in latter stages of the competition for containing challenging switching situations. Table \ref{t:networks} summarizes the scenarios considered. The two networks had 617 and 73 buses respectively, and the scenarios largely varied in terms of branches and thermal limits (all devices parameters were the same). 

\begin{table}[]
\centering
\begin{tabular}{cccccc}
\toprule
\textbf{Network} & \textbf{Scen} & \textbf{Branches} & \textbf{Devices} &  $\hat{\mathbf{P_e^{max}}}$ & $\hat{\mathbf{P_g^{max}}}$\\ \midrule
C3S4N00617D1 & 941 & 975 & 2492 & 2.80 & 0.60\\
C3S4N00617D1 & 942 & 1097 & 2492 & 2.58 & 0.60\\
C3S4N00617D1 & 943 & 1198 & 2492 & 2.43 & 0.60\\
C3S4N00617D1 & 951 & 975 & 2492 & 2.84 & 0.60\\
C3S4N00617D1 & 952 & 1097 & 2492 & 2.64 & 0.60\\
C3S4N00617D1 & 953 & 1198 & 2492 & 2.51 & 0.60\\
C3S4N00617D1 & 961 & 975 & 2492 & 2.88 & 0.60\\
C3S4N00617D1 & 962 & 1097 & 2492 & 2.71 & 0.60\\
C3S4N00617D1 & 963 & 1198 & 2492 & 2.60 & 0.60\\
C3S4N00073D2 & 991 & 129 & 571 & 3.84 & 0.24 \\
C3S4N00073D2 & 992 & 139 & 571 & 3.88 & 0.24\\
C3S4N00073D2 & 996 & 129 & 571 & 3.85 & 0.24\\
C3S4N00073D2 & 997 & 139 & 571 & 3.91 & 0.24 \\
\bottomrule   
\end{tabular}
\caption{A summary of the networks considered from the ARPA-E Grid Optimization Competition~\cite{go_networks}.}
\label{t:networks}
\end{table}

Note that the OTS problem studied here presents a subset of the one considered in the ARPA-E competition. Crucially, we consider only a single time-step, meaning we plan for the steady-state position of the network at a single instance in time. We also remove ramping limits (only applicable to multi-time-step), reserves, tap changes, and the use of capacitor banks. These later discrete actions lend themselves more readily to continuous relaxations, and thus we chose to focus on branch switching. 

\section{Results} \label{sec:res}

In this section, the results of using MGA to generate robust integer solutions to mixed integer problems are presented. We chose to focus on the OTS problem, given the availability of a number of networks which the DC-OTS problem does not generate a feasible result. However, we also present some analysis for the UC problem.

\subsection{DC Optimal Transmission Switching}

In the case of transmission switching of large networks, there will likely be multiple solutions with the optimal objective value. To see this, consider that the objective only includes generation and load costs, and that often connecting or disconnecting a single line will not interfere with generator output. Given this nature of the OTS problem, we implemented $\delta_f=0$.


We tested each MGA algorithm on each of the 13 networks previously described. For all scenarios, the DC-OTS chosen topology had no feasible AC solution. This is especially significant given that in our simulation we allow the maximum changes to generators in the AC recovery; meaning the is no generator/load configuration that results in safe voltages and branch loads.

\begin{figure}
    \centering
    \includegraphics[width=\columnwidth,clip,trim={2cm 3cm 4.7cm 3.5cm}]{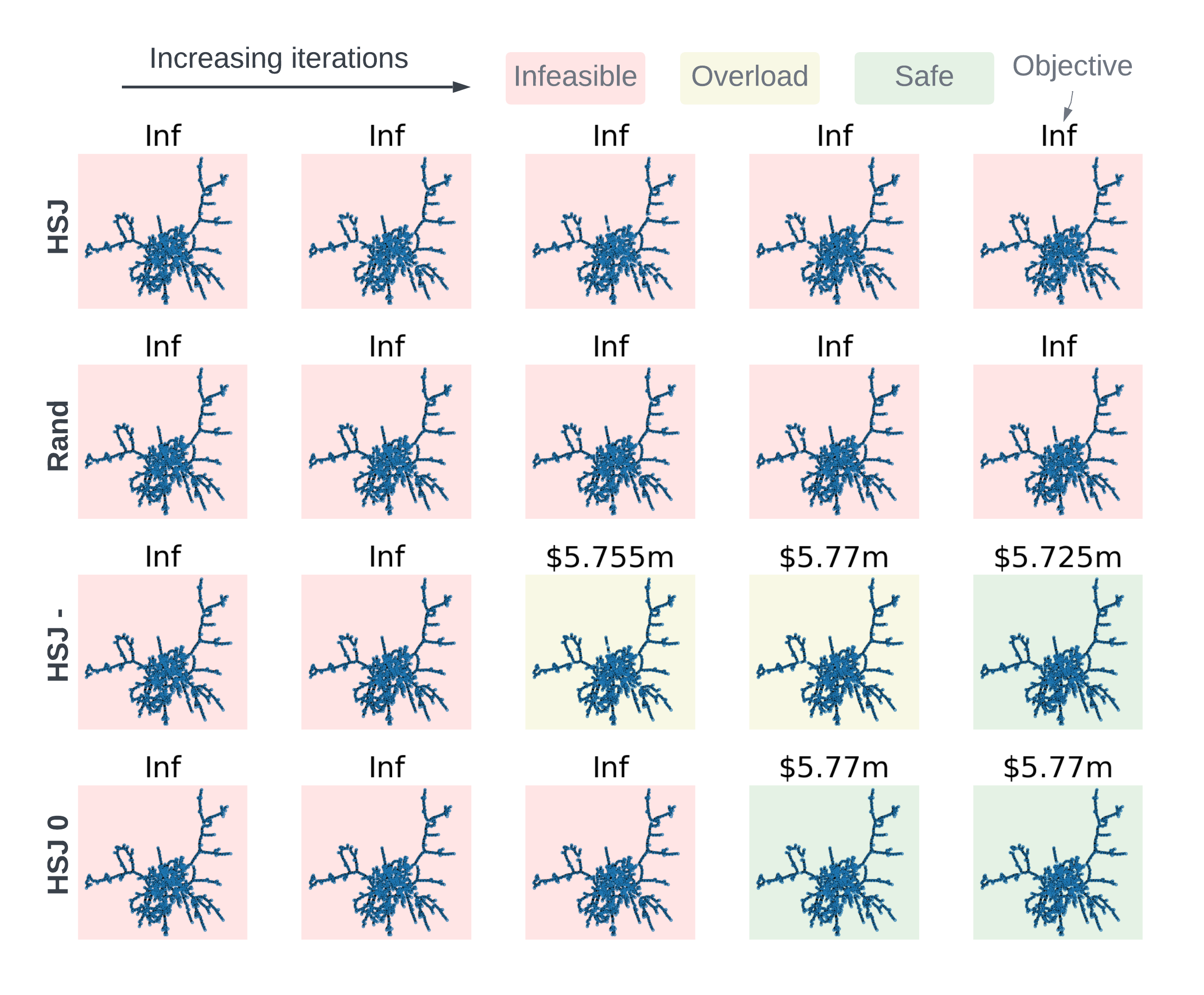}
    \caption{The chosen topologies for the first five iterations using each MGA method. The background color demonstrates whether the solution was: infeasible, had small unavoidable overloads, or was safe after an AC OPF.}
    \vspace{-.3cm}
    \label{fig:topo_demo}
\end{figure}

Figure \ref{fig:topo_demo} demonstrates the first five topologies found using each MGA method for a single scenario (note that all use the DC-OTS solution in the first iteration). We categorize these topologies according to: infeasible (no viable solution), overloads (where there are solutions which satisfy the voltage bounds but with unavoidable small branch overloads), and safe (all constraints respected). We see that the traditional HSJ and random methods do not find any feasible solution within the first 5 iterations. Meanwhile, the negative HSJ and zero sum HSJ both find feasible solutions. 

\begin{figure}
    \centering
    \includegraphics[width=\columnwidth,clip,trim={0.2cm 0.45cm 0.3cm 0.3cm}]{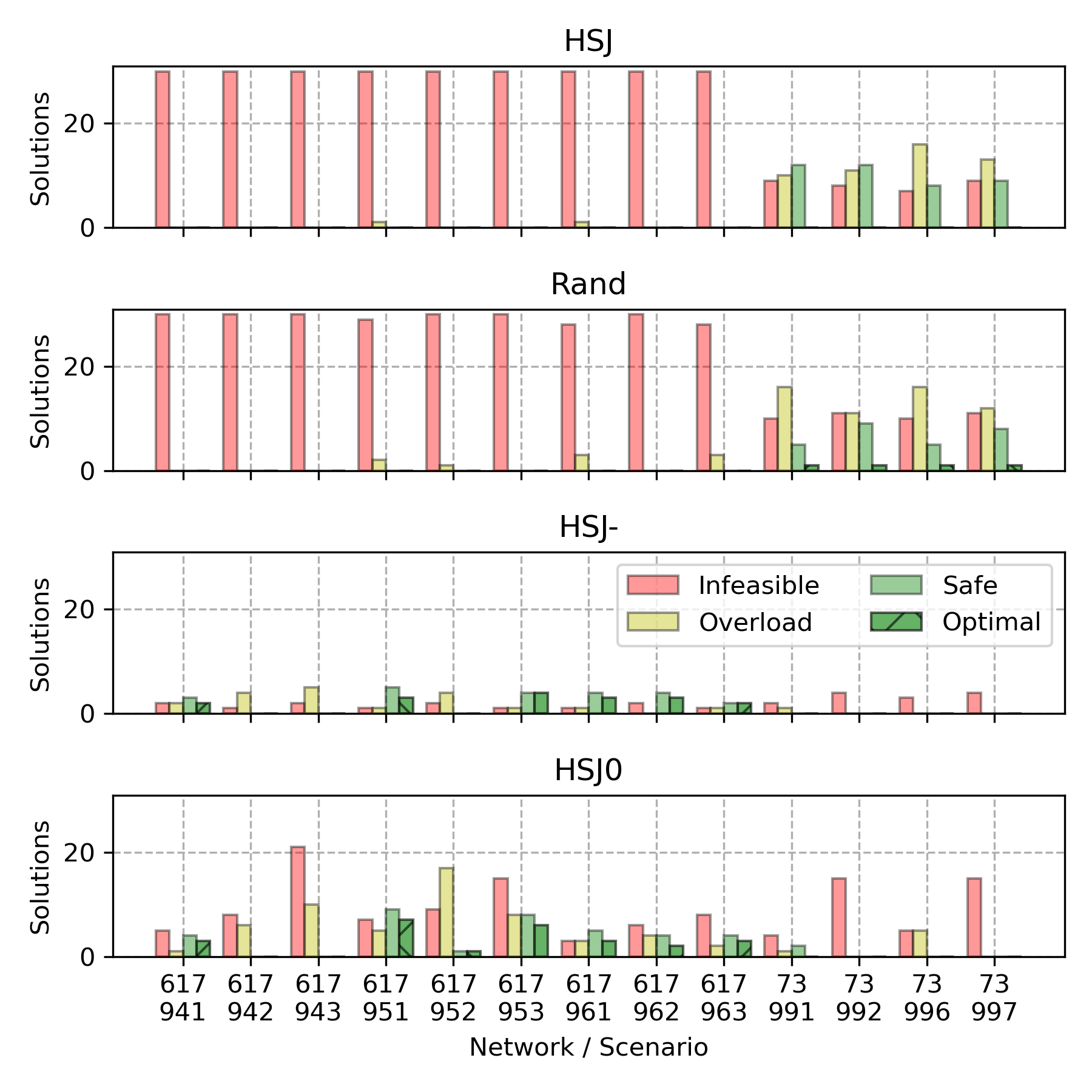}
    \caption{The breakdown of solutions into AC infeasible, overloaded, safe, and optimal after (up to) 30 iterations of each MGA method. The results are shown for each of the nine scenarios on the 617 bus network.}
    \label{fig:feas_bar}
\end{figure}

Figure \ref{fig:feas_bar} summarizes for each scenario and MGA method the break-down of the linearized solutions found in up to 30 iterations of MGA. The red bar indicates infeasible solutions -- topologies which do not have an AC power flow solution that can meet the voltage bounds. The yellow bar shows overload solutions -- which have meet the voltage bounds but suffer some branch overload. The plain green bar shows solutions which respect all constraints, and the striped green bar shows the subset of these which achieve the optimal cost.

We observe significantly different trends between the two networks. For the 617 bus network, the HSJ and random vector MGA methods generate almost exclusively infeasible solutions. The negative HSJ converges quickly (thus not generating many options) but finds feasible solutions in all cases. The zero-sum HSJ likewise always finds a feasible solution, but finds a larger number of solutions compared to the negative HSJ. On the other hand, for the 73 bus network traditional HSJ and random vectors performed best -- each identifying several safe solutions (although only random vector finds the optimal solutions). While both the negative and zero-sum methods converge quickly, mostly without finding a feasible solution.

\begin{figure}
    \centering
    \includegraphics[width=\columnwidth,clip,trim={0.2cm 0.45cm 0.3cm 0.3cm}]{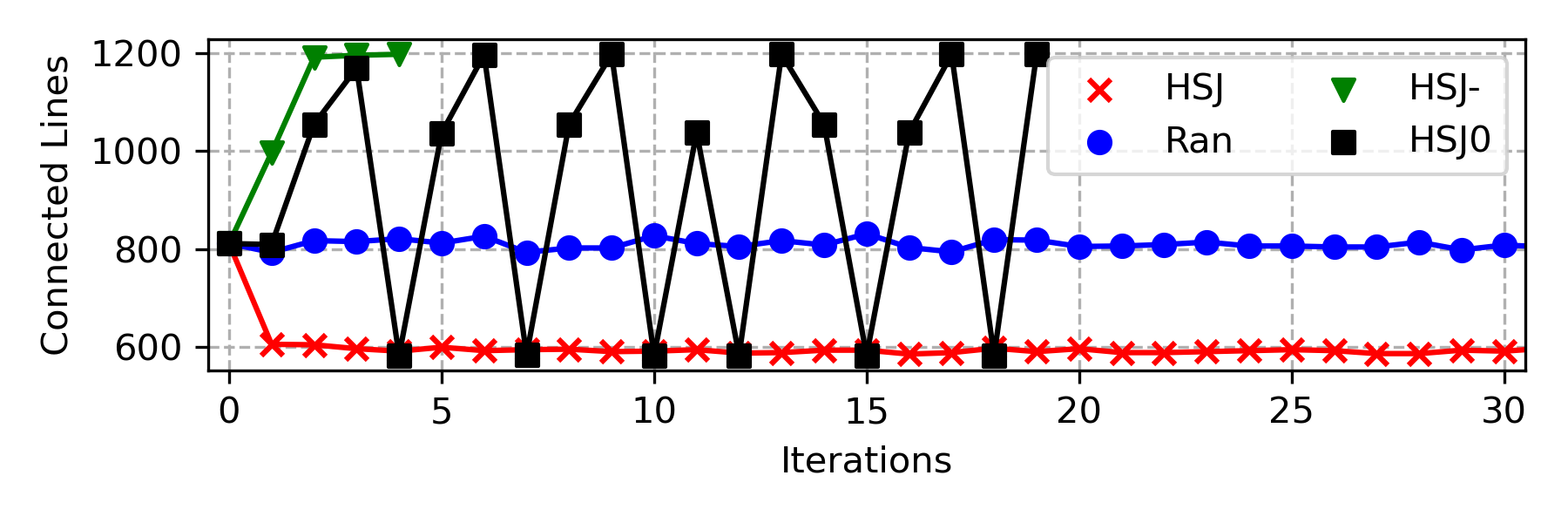}
    \caption{How the number of connected (on) branches varies between iterations of each MGA method on network 617 scenario 963.}
    \label{fig:connected_lines}
    \vspace{-.2cm}
\end{figure}

\begin{figure*}[]
    \centering
    \includegraphics[width=15cm,clip,trim={1cm 0.6cm 1cm 0.3cm}]{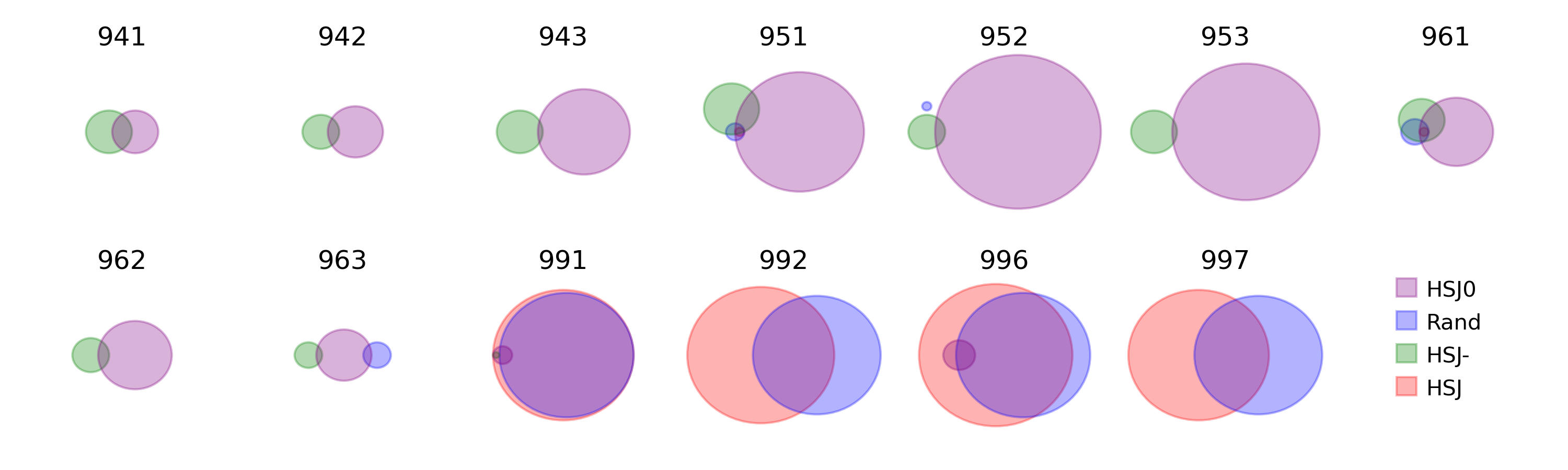}
    \caption{An illustration of the overlap between the feasible solutions determined by each method for each scenario. Note that with circles it is not always possible to display overlap exactly, in these cases the nearest approximation is used.}
    \label{fig:unique}
\end{figure*}

To better understand the differences between the solutions found be each method, Fig. \ref{fig:connected_lines} shows how the number of connected branches varies with each MGA iteration. Note that this is for a single scenario, but this was representative of the pattern seen across all scenarios and both networks. We see that (as expected) the HSJ algorithm tends to reduce the connectivity of the network with each iteration, while the negative version shows the opposite trend. The random vector method did not tend to significantly change the number of connected branches, possibly a result of the random number generation strategy used which was zero bias. The zero-sum HSJ method is the only one which explores both more and less connected strategies. From this we can infer that the 617 bus networks benefits from increased connectivity compared to the DC-OTS starting point, while the 73 bus seems to benefit from dis-connecting branches from the DC-OTS solution.

Further insight can be gained from looking at the number of feasible solutions that were found using multiple methods. If many of the same solutions are found by each method this may indicate that we are exploring the feasible space well. Figure \ref{fig:unique} visualizes the overlap between the feasible solutions generated for each scenario. The size of the circles is linearly proportional to the number of solutions found, and the amount of overlap approximately represents the number of solutions which occupy both sets. For example, for scenario 991 (of the 73 bus network) we see that HSJ manages to find almost all solutions found by other methods, with all but the random vector method being a small sub-set.

\begin{figure}
    \centering
    \includegraphics[width=\columnwidth,clip,trim={0cm 0cm 0cm 0cm}]{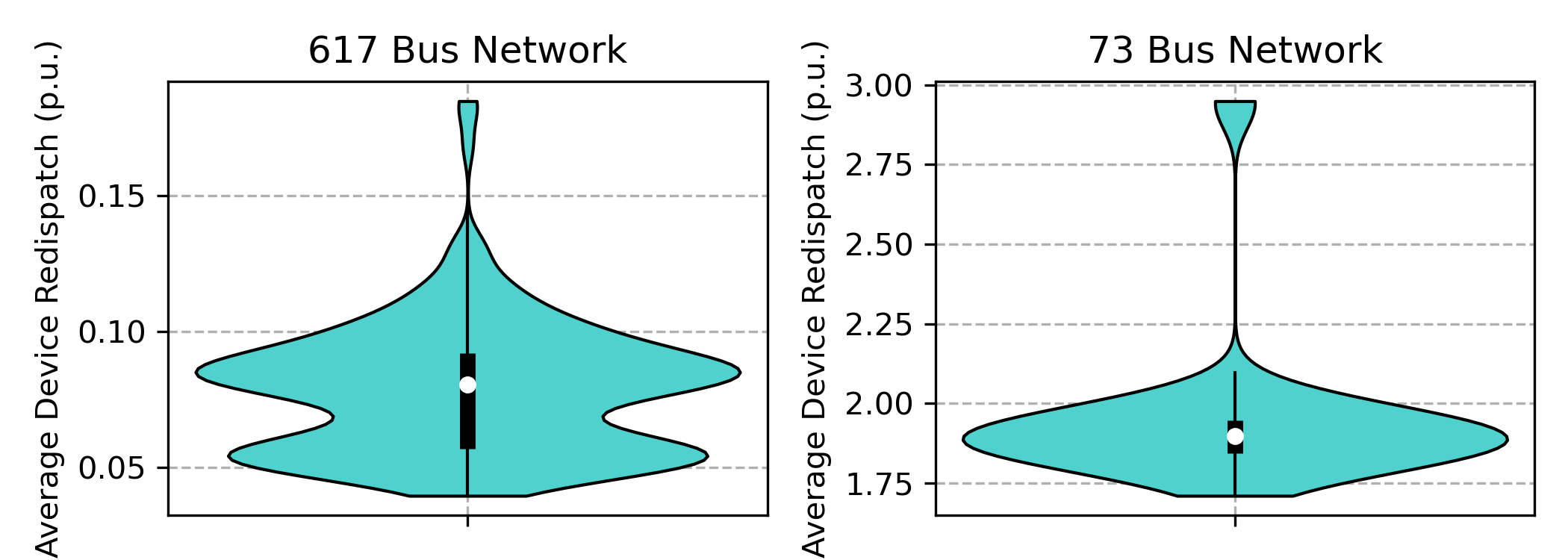}
    \caption{The distribution of the average absolute device re-dispatch between the DC-OTS and the feasible or nearest-feasible AC-OTS solution. Note that the plots have significantly different y-axis.}
    \label{fig:violin}
    \vspace{-.2cm}
\end{figure}

Again, we see a very different trend between the two networks. For the 73 bus scenarios (991--997) there is a large overlap between the solutions recovered by HSJ and random vector (and the other methods solutions tend to be a subset). Whereas, in the 617 bus scenarios we see that most solutions are only found by a single method, with small (often singular) overlaps between the zero sum and negative HSJ methods. One possible explanation for this is that the feasible region for the 73 bus network is smaller (given it has fewer branches) and thus the MGA methods have an easier time exploring. However, a more notable factor is the drastically higher percentage of solutions which turned out to be infeasible for the 617 network. 
To better understand this, we consider the amount of device (generator or load) re-dispatch between the DC-OTS and the recovery AC-OPF solutions. Figure \ref{fig:violin} shows the distribution of mean absolute device re-dispatch across all methods and scenarios for each network (note that the vertical axis has separate scales). This includes the infeasible solutions -- in this case we consider the re-dispatch that achieves the smallest sum of violations. It can be seen that, for the 617 bus network, the average re-dispatch is a factor of 10 smaller than the 73 bus network. 

\begin{figure}
    \centering
    \includegraphics[width=\columnwidth,clip,trim={0cm 0cm 0cm 0cm}]{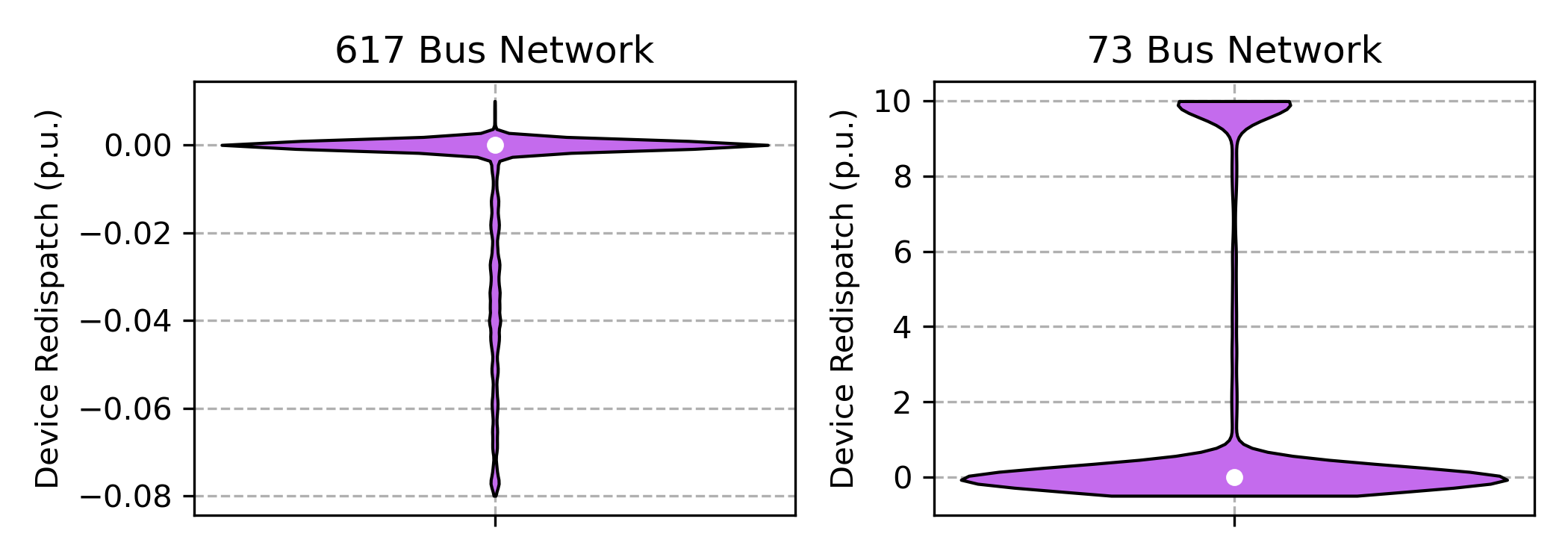}
    \caption{The distribution of individual device re-dispatch. Note that the plots have significantly different y-axis.}
    \label{fig:violin2}
    \vspace{-.1cm}
\end{figure}

Further information can be gained by looking at the changes to the individual devices, Fig. \ref{fig:violin2} shows the distribution of individual device dispatch including direction. This shows that the 617 bus network there is a net removal of both load and generation during re-dispatch, and that most devices only re-dispatch a small amount (relative to the average maximum value). Whereas, in the 73 bus network there are a small number of devices with a very high increase in output (around 10 p.u.). This demonstrates that the 73 bus network has more flexibility than the 617, due to the electrical network properties and/or device parameters.


The main benefit of the MGA solution exploration approach over previous work is that it combines the speed of linearized approaches, with the robustness of more exhaustive searching schemes. Table \ref{t:bigcomp} summarizes the number of solutions, and average computation time for each stage shown for each method and scenario. It is notable that a minimum of eight feasible AC-OTS solutions were found for each scenario, while none of these networks had a feasible topology using DC-OTS. Each iteration of the MGA and OPF is completed in seconds, allowing many solutions to explored in near-real time. For both networks the traditional HSJ method was significantly slower to iterate MGA compared to other methods, although still within reasonable limits. The average AC-OPF iterations were slower for the methods with fewer feasible solutions, likely because the problem becomes less well conditioned once the slack variables are non-zero.

\begin{table*}[t]
\centering
\begin{tabular}{cc|cccc|cccc|cccc}
\toprule
\multirow{2}{*}{\textbf{Network}} & \multirow{2}{*}{\textbf{Scenario}} & \multicolumn{4}{c}{\textbf{\# Feasible Solutions}} & \multicolumn{4}{c}{\textbf{Average MGA Iteration Time (s)}} & \multicolumn{4}{c}{\textbf{Average OPF Iteration Time (s)}} \\
                                  &                                    & HSJ        & HSJ-        & HSJ0       & Rand       & HSJ         & HSJ-         & HSJ0         & Rand        & HSJ         & HSJ-         & HSJ0         & Rand        \\ \midrule
         C3S4N00617D1 & 941 & 0& 5& 5& 0& 2.07 & 0.29 & 0.34 & 0.29 & 2.14 & 1.20 & 1.39 & 1.40 \\
         C3S4N00617D1 & 942 & 0& 4& 6& 0& 1.76 & 0.27 & 0.36 & 0.30 & 2.26 & 1.29 & 1.51 & 1.61 \\
         C3S4N00617D1 & 943 & 0& 5& 10& 0& 1.96 & 0.41 & 0.46 & 0.38 & 2.09 & 1.34 & 1.47 & 1.70 \\
         C3S4N00617D1 & 951 & 1& 6& 14& 2& 2.10 & 0.31 & 0.43 & 0.28 & 2.09 & 1.19 & 1.56 & 1.43 \\
         C3S4N00617D1 & 952 & 0& 4& 18& 1& 1.78 & 0.40 & 0.40 & 0.32 & 2.43 & 1.28 & 1.45 & 1.59 \\
         C3S4N00617D1 & 953 & 0& 5& 16& 0& 2.01 & 0.28 & 0.39 & 0.35 & 2.08 & 1.33 & 1.48 & 1.69 \\
         C3S4N00617D1 & 961 & 1& 5& 8& 3& 2.07 & 0.26 & 0.30 & 0.27 & 2.15 & 1.18 & 1.40 & 1.37 \\
         C3S4N00617D1 & 962 & 0& 4& 8& 0& 1.92 & 0.27 & 0.34 & 0.31 & 2.34 & 1.33 & 1.40 & 1.55 \\
         C3S4N00617D1 & 963 & 0& 3& 6& 3& 1.89 & 0.27 & 0.40 & 0.35 & 2.13 & 1.36 & 1.54 & 1.76 \\ \midrule
         C3S4N00073D2 & 991 & 22& 1& 3& 21& 0.08 & 0.04 & 0.05 & 0.05 & 0.34 & 0.27 & 0.25 & 0.24 \\
         C3S4N00073D2 & 992 & 23& 0& 0& 20& 0.11 & 0.05 & 0.05 & 0.05 & 0.21 & 0.25 & 0.26 & 0.24\\
         C3S4N00073D2 & 996 & 24& 0& 5& 21& 0.08 & 0.04 & 0.05 & 0.05 & 0.23 & 0.27 & 0.23 & 0.24 \\
         C3S4N00073D2 & 997 & 22& 0& 0& 20& 0.11 & 0.05 & 0.05 & 0.05 & 0.23 & 0.26 & 0.26 & 0.25\\
         \bottomrule
\end{tabular}
\caption{For each of the tested networks: number of feasible solutions found before convergence or the 30th iteration (whichever was soonest), average time to compute each MGA solution, and average time to validate the solution in AC.}
\label{t:bigcomp}
\end{table*}

\subsubsection{Comparison against heuristic methods}

\begin{figure}
    \centering
    \includegraphics[width=\columnwidth,clip,trim={0cm .4cm 0cm 0cm}]{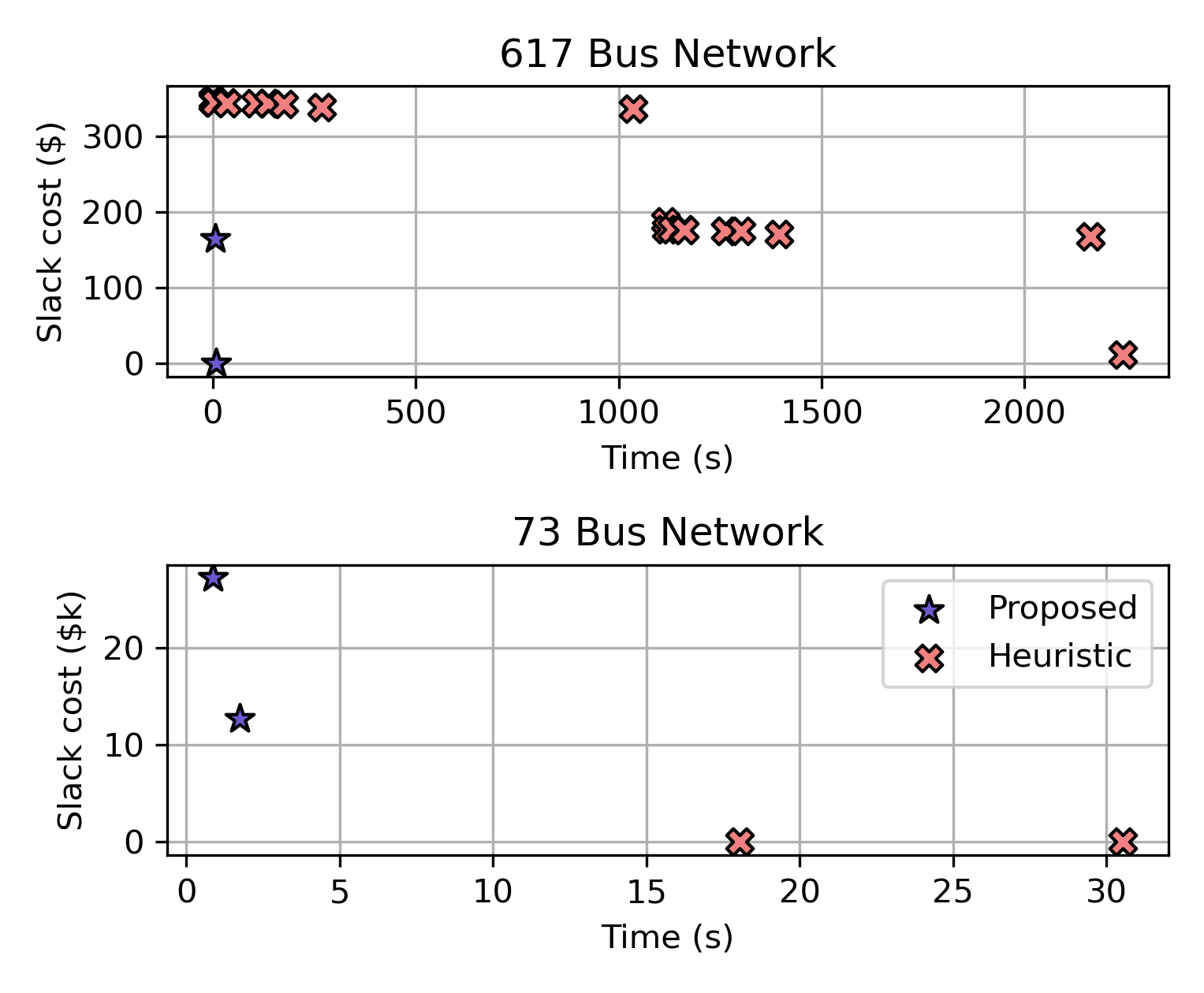}
    \caption{A comparison of the convergence time of the proposed method to a greedy heuristic algorithm, where in each iteration the single best line switch is found through exhaustive search.}
    \label{fig:heuristic_comp}
    \vspace{-.2cm}
\end{figure}

Here we compare the proposed method against a popular greedy heuristic algorithm; where the best single topology switch is found by exhaustion, executed, and then search is repeated. This algorithm converges when no switch yields an improvement, or the improvement of the best switch is below a pre-defined tolerance. This algorithm can not guarantee the global optima, but does guarantee AC feasibility and that the final topology is no worse than the initial topology. Figure \ref{fig:heuristic_comp} visualizes how the optimality gap (measured in dollars of sub-optimality) decreases clock-time for the greedy algorithm compared to the proposed method. 

For the 617 bus network we see that the heuristic method finds the optimal solution within seconds, while the greedy algorithm takes more than 30 minutes to find a comparable solution. This is due to the much larger number of AC power flow calculations that the greedy method needs to execute -- each iteration requires $n_e$ calculations, while in the proposed method there is only one AC power flow solution per iteration. For the 73 bus network we see that the proposed method converges within two seconds, but does not find the globally optimal solution. Whereas the greedy algorithm finds a superior solution in around 30 seconds (in this case the AC power flow calculations are less costly due to the smaller network). This demonstrates that the proposed method converges orders of magnitude faster than this existing heuristic, although better solutions may be found with more exploratory search criteria.

\subsection{DC Unit Commitment}

We believe that MGA may be useful in a variety of integer decisions where the MILP formulation may lead to infeasible solutions. Here we demonstrate the application of MGA to the DC unit commitment problem. Using scenario 991 of the aforementioned 73 bus network we used MGA to search for alternative single-timestep dispatch scenarios. Unlike in OTS, it is unlikely that many dispatch solutions with the same objective function exist, so we consider two values for $delta$, \$100 and \$10,000. The Random Vector search algorithm was chosen in order to maximally explore the solution space.

\begin{figure}
    \centering
    \includegraphics[width=\columnwidth,trim={2cm 3cm 2cm 3cm},clip]{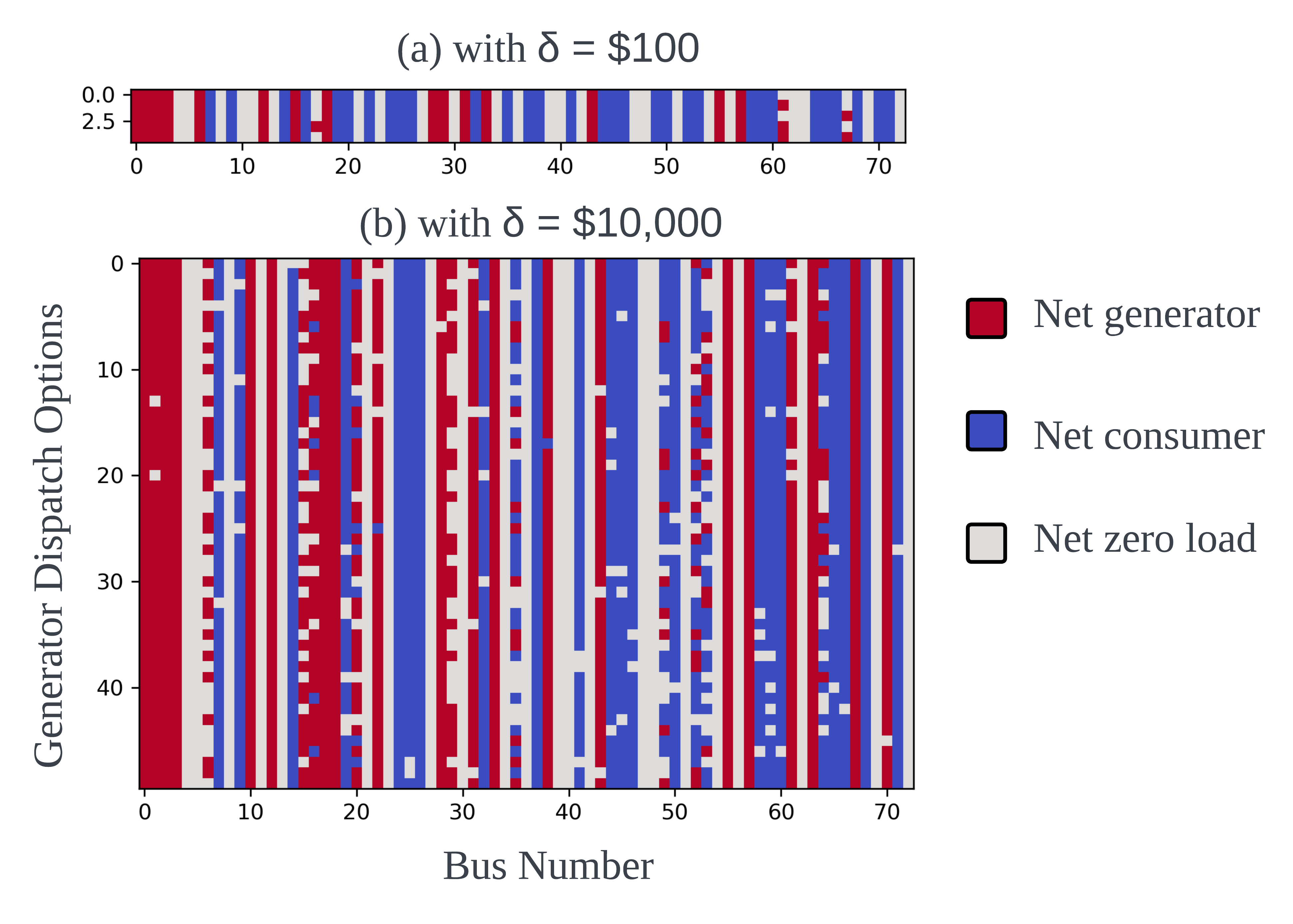}
    \caption{A visualization of the generator dispatch scenarios discovered for two different objective sensitivity values.}
    \label{fig:dcuc}
\end{figure}

Figure \ref{fig:dcuc} visualizes whether each bus becomes a net generator, load, or junction (as a reminder, we have dispatchable demand as well as supply in this case). Only unique dispatches are shown (the random vector algorithm sometimes generates repeated solutions). We see that the number of distinct solutions increases significantly for larger tolerances. For the \$100 case (equivalent to $\sim$0.002\% tolerance) only five distinct solutions were found, whereas for \$10,000 ($\sim$2\%) 50 solutions were found.

In this case, the DC-UC problem produced an AC feasible solution. However, there are cases where the thermal limits interfere with the DC-UC solution; overall 43 of the 55 generated dispatch solutions were feasible. Feasibility becomes more challenging in the multi-period problem, in which case the percentage of dispatches which are feasible is likely to be smaller.

\section{Discussion \& Conclusions} \label{sec:con}

This paper investigated the use of modeling to generate alternatives (MGA) to improve the success of linearized mixed-integer OPF problems. Specifically, we considered the optimal transmission switching (OTS) and unit commitment problems. We investigated the use of two classical MGA methods: hop, skip and jump (HSJ) and random vector, as well as proposing two novel modifications to the HSJ method. Here we considered two realistic networks that had challenging transmission switching solutions, where DC-OTS yielded an infeasible solution. Each network had multiple scenarios, where only the available branches and their limits was altered, a total of 13 scenarios were tested. 

We showed that in every scenario we were able to find several good AC-OTS solutions using one of the MGA methods. This is significant, given that in all cases the DC-OTS solution was infeasible. However, the performance between the MGA methods varied significantly between networks. One network (with 73 buses) appeared to have a large amount of flexibility, and benefited from removing lines compared with the DC-OTS solution. In this case we found the best performance with the traditional MGA methods (HSJ and random vector). The other network was highly constrained, and benefited from adding branches. In this case, the traditional methods failed to find a good solution, while the newly proposed methods (negative and zero-sum HSJ) found a variety of good solutions. Further work is needed to understand how MGA methods can be selected based on network parameters.



\bibliography{refs.bib}{}
\bibliographystyle{IEEEtran}




\end{document}